\documentstyle[12pt]{article}
\newtheorem{defin}{Definition}
\newtheorem{prop}{Proposition}
\newtheorem{nt}{Remark}
\newtheorem{Th}{Theorem}

\newtheorem{ex}{Example}

\newfont{\sdbl}{msbm9}
\newfont{\dbl}{msbm10 at 12pt}
\newcommand{\eqdef}{\stackrel{\rm def}{=}}
\newcommand{\proof}{{\bf Proof.\ }}

\newcommand{\oo}{{\cal O}}
\newcommand{\hoo}{\hat{\oo}}

\newcommand{\g}{{\cal G}}
\newcommand{\h}{{\cal H}}
\newcommand{\ad}{{\cal A}}

\newcommand{\mod}{\mathop {\rm mod}}

\newcommand{\Div}{\mathop {\rm Div}}

\newcommand{\Lim}{\mathop {\rm lim}}

\newcommand{\nm}{\mathop {\rm Nm}}
\newcommand{\Spec}{\mathop {\rm Spec}}
\newcommand{\sign}{\mathop {\rm sign}}

\newcommand{\Frac}{\mathop {\rm Frac}}
\newcommand{\dm}{\mathop {\rm dim}}

\newcommand{\da}{{\mbox{\dbl A}}}
\newcommand{\dz}{{\mbox{\dbl Z}}}
\newcommand{\dbz}{\dz}

\newcommand{\dbq}{{\mbox{\dbl Q}}}
\newcommand{\dbf}{{\mbox{\dbl F}}}

\newcommand{\f}{{\cal F}}
\newcommand{\lto}{\longrightarrow}

\newcommand{\dn}{{\mbox{\dbl N}}}
\newcommand{\sdbn}{{\mbox{\sdbl N}}}
\newcommand{\df}{{\mbox{\dbl F}}}
\newcommand{\D}{{\bf D}}
\newcommand{\C}{{\bf C}}

\newcommand{\lm}{\mathop{\rm lim}}

\newcommand{\pl}{\lm\limits_{\longleftarrow}}
\newcommand{\QS}{{\rm \bf QS}}
\newcommand{\Sh}{{\rm \bf Sh}}
\newcommand{\Ab}{{\rm \bf Ab}}

\begin{document}
\author{D.V. Osipov \footnote{The author was supported by an LMS grant for young
Russian mathematicians at the University of Manchester, and also by the Russian Foundation
for Basic Research, grant no. 05-01-00455.}}

\title{$n$-dimensional local fields and adeles on $n$-dimensional schemes.}
\date{}
\maketitle
\section{Introduction}

The notion of $n$-dimensional local field  has appeared in the works of
A.~N.~Parshin and K.~Kato  in the middle of 70's.
These fields generalize the usual local fields (which is $1$-dimensional
in this sense) and help us to see on higherdimensional
algebraic schemes from the local point of view.

With  every flag
$$ X_0 \subset X_1 \ldots \subset
X_{n-1}$$
$$
{\rm dim \,} X_i = i $$
of irreducible subvarieties on a scheme $X$
($
{\rm dim \,} X =n
$)
one can canonically associate
a ring $K_{(X_0, \ldots, X_{n-1})}$.
  In case  everything is
regularly embedded, the ring is
an $n$-dimensional local field.

Originally, higherdimensional local fields
were used for the  development of generalisation
of class field theory to the schemes
of arbitrary dimension (works of A.~N.~Parshin,
K.~Kato, S.~V.~Vostokov and others), \cite{P}, \cite{KS}.
But many problems of  algebraic varieties
can be reformulated in  terms
of higherdimensional local fields
and higher adelic theory.

For a scheme $X$ adelic object
$$ \da_X =
\prod\nolimits' K_{(X_0, \ldots, X_{n-1})}
$$
where the product is taken over all the flags
with respect to certain restrictions on components of adeles.
A.N. Parshin defined in~\cite{P1} adeles on algebraic surfaces, which generalize  usual adeles on curves. A.A. Beilinson introduced a simplicial approach to adeles and
generalized to arbitrary dimensional Noetherian schemes in~\cite{B}.

A.~N.~Parshin, A.~A.Beilinson, A.~Huber,
A.~Yekutiely, V.~G.~Lomadze  and others  described
connections of higher adelic groups with  cohomology of
coherent sheaves (\cite{P1}, \cite{B}, \cite{H}, \cite{Y}, \cite{F}, \cite{FP}),
intersection theory (\cite{P2}, \cite{Lom}, \cite{Os0}, \cite{FP}),
Chern classes (\cite{P2}, \cite{HY}, \cite{FP}), theory of residues (\cite{P1}, \cite{Y},
\cite{B}, \cite{L}, \cite{FP}), torus actions (\cite{GP}).

This paper is a survey of basic notions of higherdimensional local fields
and adeles on higherdimensional schemes.

The paper is organized as follows.

In section \ref{sect2}
we give a general definition of $n$-dimensional local field
and formulate classification theorems of $n$-dimensional local fields.
We describe how $n$-dimensional local fields appear from algebraic varieties
and arithmetical schemes.

In section \ref{adel}
we define higher dimensional adeles and adelic complexes.
Starting from an example of adelic complexes on algebraic curves, we give
a general simplicial definition for arbitrary Noetherian schemes, which is due
to A.A. Beilinson. We formulate the theorems about adelic resolutions of quasicoherent sheaves
on Noetherian schemes. We applicate this general constructions to algebraic sufaces
to obtain adelic complexes on algebraic surfaces, which were introduced by A.N.~Parshin.

In section \ref{sect4} we describe restricted adelic complexes.
In constrast to the adelic complexes from section~\ref{adel}, restricted adelic complexes are connected with a single flag of subvarieties. A.N. Parshin introduced restricted adeles
for algebraic surfaces in~\cite{P5}, \cite{P4}. The author introduced restricted adelic
complexes for arbitrary schemes in \cite{Os}. We give  also the reconstruction theorem on
restricted adelic complexes.

In the last section we briefly describe reciprocity laws on algebraic surfaces.

The author  is very grateful to A.N.~Parshin for a lot of discussions on higherdimensional local fields and adeles. The author is also grateful to M.~Taylor for interesting discussions and his hospitality during the visit to Manchester sponsored by an LMS grant.

\section{$n$-dimensional local fields} \label{sect2}
\subsection{Classification theorems.}

We fix a perfect field $k$.

We say that $K$ is a local field of dimension $1$ with the residue field $k$,
if $K$ is a fraction field of complete discrete valution ring $\oo_K$
with the residue field $\bar K =k$.

We denote by $\nu_K$ the discrete valuation of $K$
and  by $m_K$ the maximal ideal of ring $\oo_K$.

Such a field has the following structure

$$ K  \supset \oo_K \to \bar K = k$$.

As examples of such fields we have the field of power series
$$ K = k((t)) \mbox{,}    \qquad  \oo_K = k [[t]] \mbox{,}  \qquad \bar K = k $$
and the field of p-adic numbers
$$ K = \dbq_p  \mbox{,}   \qquad \oo_K = \dbz_p  \mbox{,}  \qquad  k = \dbf_p  \mbox{.}$$

Moreover, we have only the following possibilities,
see~\cite[ch. II]{Ser}:
\begin{Th}  \label{cl1}
Let $K$ be a local field of dimension $1$ with the residue field $k$ then
\begin{enumerate}
    \item  \label{c1}
     $K= k((t))$ is the  power series field if $char K = char k$;
    \item  \label{c2}
    \begin{enumerate}
    \item
    $K = \Frac (W(k))$ where $\oo_K = W(k)$ is the Witt ring
(for example, $K = \dbq_p$),
         \item or $K $ is  a finite totally
ramified extension of the field $\Frac (W(k))$
         \end{enumerate}
     if $char K = 0$, $char k =p$.
\end{enumerate}
\end{Th}

Now we give the following inductive definition.

\begin{defin}
We say that a field $K$ is a local field of dimension $n$
with the last residue filed $k$ if
\begin{enumerate}
\item $n = 0$ and $K =k$
\item $n \ge 1$ and $K$ is the fraction field of a complete discrete valuation
ring $\oo_K$ whose residue field $\bar K$ is a local field of
dimension $n-1$ with  the last residue field $k$
\end{enumerate}
\end{defin}

A local field of dimension $n \ge 1$ has the following inductive structure:
$$ K = K^{(0)} \supset \oo_K \to \bar{K} = K^{(1)} \supset \oo_{\bar{K}}  \to
\bar{K}^{(1)}= K^{(2)} \supset \oo_{K^{(2)}} \to \ldots \bar{K}^{(n)}=k  \mbox{,} $$
where for a discrete valuation field $F$ $\oo_F$  the ring
of integers in $F$ and $\bar{F}$ the residue field.
The maximal ideals in $\oo_{K^{(i)}}$ we denote by $m_{K^{(i)}}$.

Every field $K^{(i)}$ is a local field of dimension $n-i$
with the last residue field $k$.

\begin{defin}
A collection of elements $t_1, \ldots, t_n \in \oo_K$
is called  a system of local parameters, if for all $i = 1, \ldots, n$
$$
t_i \; \mod m_{K^{(0)}} \in \oo_{K^{(1)}}, \; \ldots, \;
t_i \; \mod m_{K^{(i-2)}}
\in \oo_{K^{(i-1)}}
$$
and the  element
$t_i \: \mod m_{K^{(i-2)}}
\in \oo_{K^{(i-1)}}$
is a generator of $m_{K^{(i-1)}}$.
\end{defin}

An example of an $n$-dimensional local field is the field $K = k ((t_n)) \ldots ((t_1))$.
For this field we have
$$
K^{(i)}=k((t_n)) \ldots ((t_{i+1}))
\qquad
\mbox{,}
\qquad
\oo_{K^{(i)}}= k ((t_n)) \ldots ((t_{i+2}))[[t_{i+1}]] \mbox{.}
$$

We consider $n =2$, then we have
$$K \supset \oo_K \to \bar{K} \supset \oo_{\bar{K}} \to k \mbox{.}$$

We can construct the following examples of $2$-dimensional local fields with the last residue field $k$. These examples depend on  characteristic of $2$-dimensional field $K$
and  characteristics of its residue fields.
\begin{itemize}
\item $K = k ((t_2))((t_1))$. And $char K = char \bar{K} = char k$.
\item $K = F ((t))$, where $F$ is a local field of dimension $1$ with the residue field $k$
such that $char(F) \ne char(k)$, for example $ F = \dbq_p$.
\item $K = F \{\{t\}\}$, where $F$ is a local field of dimension $1$ with the residue field $k$.
\end{itemize}

The field $F \{\{t\}\}$ has the following description
$$
a \in F\{\{t\}\} \quad : \quad a= \sum_{i= - \infty}^{i = + \infty} a_i t^i
\mbox{,} \qquad a_i \in F \mbox{,}
$$
$$
\mbox{where} \qquad         \mathop{\lim}\limits_{i \to -\infty} \nu_F(a_i) = 0 \qquad \mbox{and}
\qquad          \nu_F(a_i) > c_a \qquad  \mbox{for some integer} \qquad c_a \mbox{.} $$

We put the discrete valuation $\nu_{F\{\{t\}\}} (a) = \min \nu_{F} (a_i)$.
Then the ring
$ \oo_{F\{\{t\}\}}$ consits of elements $a$ such that all $a_i \in \oo_F $, and
the maximal ideal
$ m_{F\{\{t\}\}}$   consists of elements $a$ such that all $ a_i \in m_F                                       $.
Therefore
$$\overline{F\{\{t\}\}} = \bar{F}((t)) \mbox{.}$$

We remark that for $F = k((u))$
the field  $F\{\{t\}\}$ is isomorphic to the field  $k ((t))((u))$.

There is the following classification theorem, see \cite{FP}, \cite{P3}, \cite{Zh}, \cite{Zh1}.

\begin{Th} \label{cl2}
Let $K$ be an $n$-dimensional local field with the finite last residue field.
Then
\begin{enumerate}
\item
   \label{ca1}
 $K$ is isomorphic to $\dbf_q((t_n)) \ldots ((t_1))$ if $char (K) = p$;
\item
   \label{ca2}
$K$ is isomorphic to $F ((t_{n-1})) \ldots ((t_{1}))$,
$F$ is a $1$-dimensional local field, if $char (K^{(n-1)}) = 0$.
\item
   \label{ca3}
$K$ is a finite extension of a  field
$$
F \{ \{ t_n \} \}  \ldots \{ \{ t_{m+2} \} \} ((t_m)) \ldots ((t_1))  \qquad   \qquad \qquad (*)
$$
and there is a finite extension of $K$ which is of the form $(*)$, but possibly
with different $F$ and $t_i$, if $char (K^{(m)}) = 0$, $char  (K^{(m+1)}) = p$.
\end{enumerate}
\end{Th}

We remark that if $\pi$ is a local parameter for a $1$-dimensional local field $F$
then $t_1, \ldots, t_m, \pi, t_{m+2}, \ldots, t_n$
are local parameters for a field
$$  F \{ \{ t_n \} \}  \ldots \{ \{ t_{m+2} \} \} ((t_m)) \ldots ((t_1)) \mbox{.}$$

\subsection{Local fields which come from algebraic geometry.}
We consider an algebraic curve $C$ over the field $k$.
We fix a smooth point $p \in C$. We consider
$$
K_p = \Frac ( \hat{\oo}_p)
\mbox{,}$$
where
$$ \hat{\oo}_p =  \mathop{\mathop{\lim}_{\leftarrow}}_n \oo_p/ m_p^n
$$
is a completion of the local ring $\oo_p$ of the point $p$ on the curve $C$.
Then
\begin{equation} \label{f1}
K_p = k(p) ((t))  \mbox{,}
\end{equation}
where $k(p) = \oo_p/ m_p $ is the residue field of the point $p$ on the curve $C$,
which is the finite extension of the field $k$, and $t$ is a local parameter of the point
$p$ on the curve $C$.

We see that the field $K_p$ corresponds to the case~\ref{c1}
of the classification theorem~\ref{cl1}.

\vspace{0.5cm}

Now we consider a field of algebraic numbers $K$,
which is a finite extension of the field $\dbq_p$. We consider the ring of integers
$A$ of the field $K$. Let $X = \Spec A$ be a $1$-dimensional scheme.
We fix a closed point $p \in X$, which corresponds to a maximal ideal in $A$.
Then the completion of the field $K$ at the point $p$ is
\begin{equation} \label{f2}
K_p = \Frac \,( \, \mathop{\mathop{\lim}_{\leftarrow}}_n A_p/ m_p^n \,)
\end{equation}

We see that $K_p$ is a $1$-dimensional local field with the residue field $\dbf_q$,
and the field $K_p$ corresponds to the case~\ref{c2}
of the classification theorem~\ref{cl1}.

\vspace{1cm}

We give now the definitions for the general situation.

Let $R$ be a ring, $p$ a prime ideal of $R$, $M$ an $R$-module.

Let $S_p = R \setminus p$. We write $S_p^{-1} M$ for the localisation of $M$ at $S_p$.

For an ideal $a$ of $R$ set
$C_a M =
\mathop{\mathop{\lim}\limits_{\leftarrow}}\limits_{n \in \sdbn} M / a^n M
$.

Let $X$ be a Noetherian scheme of dimension $n$.
Let $\delta = (p_0, \ldots, p_n)$ be a chain of points of $X$ (i.e. the chain of integral irreducible subschemes if we consider the closures of points $p_i$) such that $p_{i+1} \in \overline{\{p_i\}}$ for any $i$ ($ \overline{\{p_i\}} $ is the closure of the point $p_i$ in $X$).
We suppose that for all $i$ $\dim p_i = i$.
We restrict $\delta$ to some affine open neighbourhood $\Spec B$ of the closed point
$p_n$ on $X$.  Then $\delta$ determines a chain of prime divisors of the ring $B$,
which we denote by the same letters $(p_0, \ldots, p_n)$.
We define a ring
\begin{defin}
\begin{equation} \label{f}
K_{\delta} \eqdef C_{p_0}S_{p_0}^{-1} \ldots C_{p_n} S_{p_n}^{-1} B
\end{equation}
\end{defin}

This definition of $K_{\delta}$ does not depend on the choice
of affine neighbourhood $\Spec B$ of the point $p_n$ on  the scheme $X$,
see~\cite[prop.3.1.3., prop. 3.2.1.]{H}.

We remark that the ring  $C_{p_n} S_{p_n}^{-1} B$ from formula \ref{f}
coincides with the completion $\hat{\oo}_{p_n, X}$ of the local ring of the point $p_n$
on the scheme $X$.

\vspace{1cm}

We consider now examples of formula~(\ref{f})
for small $n$.

\begin{ex} \em
If $X$ is an irredusible $1$-dimensional scheme ($X$ is an irreducible curve over the field
$k$ or the spectrum of the ring of algebric integers), $p$ is a smooth point of $X$, $\eta$ is a general point of $X$, then for $\delta = (\eta, p)$ we obtain that $K_{\delta}$
is a $1$-dimensional local field, which coincides with the field $K_p$
from formula~(\ref{f1}) or (\ref{f2}).
\end{ex}

\begin{ex} \em
Now let $X$  be an irreducible algebraic surface over the field $k$.
Let $C$ be an irreducible divisor on $X$ and $p$ a point on $C$.
We suppose that $p$ is a smooth point on $X$ and on $C$. Let $\eta$
be a general point of $X$. We consider $\delta = (\eta, C, p)$.
We fix the local parameter $t \in k(X)$ of the curve $C$ on $X$ at the point $p$ ($t = 0$ is a local equation of the curve $C$ at the point $p$ on $X$),
and local parameter $u \in k(X)$ at the point $p$ on $X$ which is transversal to the local parameter $t$  (the divisor $u = 0$ is transversal to the divisor $t = 0$ at the point $p$).
We fix any affine neighbourhood $\Spec B$ of $p$ on $X$. Then
$$ C_{p} S_{p}^{-1} B = k(p) [[u,t]] $$
$$ C_{C}  S_{C}^{-1}  C_{p} S_{p}^{-1} B = k(p) ((u))[[t]] $$
$$ K_{\delta} = C_{\eta}  S_{\eta}^{-1}
 C_{C}  S_{C}^{-1}  C_{p} S_{p}^{-1} B = k(p) ((u)) ((t)) \mbox{.}
$$

Hence $K_{\delta}$ is a $2$-dimensional local field  with the last residue field $k(p)$.

We see that the field $K_{\delta}$
corresponds to the case~\ref{ca1} of the classification theorem~\ref{cl2}.
\end{ex}

\begin{ex} \em
The previous example can be generalized. Let $p_0, \ldots, p_n$
be a flag of irreducible subvarieties on $n$-dimensional alfebraic variety $X$ over the field $k$ such that  $\dim p_i = n-i$, $p_{i+1} \subset p_i$ for all $i$ and the point $p_n$ is a smooth point on all subvarieties $p_i$. We can choose a system of local parameters $t_1, \ldots,
 t_n \in \oo_{p_n, X}$ of the point $p_n$ on $X$ such that for every $i$
equations  $t_1 = 0, \ldots, t_i = 0 $ define a subvariety $p_i$ in some neighbourhood
of the point $p_n$ on $X$. Then according to formula~(\ref{f}) and similar to the previous example we have for $\delta = (p_0, \ldots, p_n)$
$$
K_{\delta} = k(p) ((t_n)) \ldots ((t_1)) \mbox{.}
$$

\end{ex}

\begin{ex} \em
Now we suppose that a scheme $X$ is an arithmetical surface,
i.e., $\dim X = 2$ and we have a flat, projective morphism
$f : X \to Y = \Spec A$, where $A$ is the ring of integers of a number field $K$.
We consider two kinds  of integral irreducible $1$-dimensional closed subschemes $C$ on $X$.
\begin{enumerate}
\item
A subscheme $C$ is horizontal, i.e., $f(C) = Y$.
We consider a point
$x \in C$ which is  smooth  on $X$ and $C$. Let $\delta = (\eta, C, x)$,
where $\eta$ is a general point of $X$. Then
$$
K_{\delta} = L((t)) \mbox{,}
$$
where $t= 0$ is a local equation of $C$ at the point $x$ on $X$ and
$L \supset K_{f(x)} \supset \dbq_p$ is a finite extension.
Thus, $K_{\delta}$
is a $2$-dimensional local field with the finite last residue field.

We see that this field $K_{\delta}$ corresponds to the case~\ref{ca2}
of the classification theorem~\ref{cl2}.
\item
A subscheme $C$ is vertical, i.e. it is a component of a fibre of $f$.
This $C$ is defined over some finite field $\df_q$.
We consider a point
$x \in C$ such that the  morphism $f$ is smooth at $x$
and the point $x$ is also defined over the field $\df_q$. Let $\delta = (\eta, C, x)$,
where $\eta$ is a general point of $X$.
Then we applicate  formula (\ref{f}). For any affine neighbourhood $\Spec B$ of $x$ on $X$
 the ring
$
 C_{x} S_{x}^{-1} B $
coincides with the completion $\hat{\oo}_{x, X}$ of the local ring of the point $x$ on $X$.
But since $f$ is a smooth map at $p$, the ring
$$
 \hat{\oo}_{x, X}
= \oo_{K_{f(x)}}[[u]]  \mbox{.}
$$
Therefore we obtain
$$
K_{\delta} = K_{f(x)}\{\{u\}\} \mbox{,}
$$
where $K_{f(x)} \supset \dbq_p $
is a finite extension.
Thus, $K_{\delta}$
is a $2$-dimensional local field with the finite last residue field.

We see that this field $K_{\delta}$ corresponds to the case~\ref{ca3}
of the classification theorem~\ref{cl2}.
\end{enumerate}

We remark that in these both cases we have a canonical embedding $f^*$ of the $1$-dimensional local field $K_{f(x)}$
to the $2$-dimensional local field
$K_{\delta}$.
\end{ex}

\vspace{0.5cm}

Now we consider only excellent Noetherian schemes $X$
(e.g. a scheme of finite type over a field, over $\dz$,
or over a complete semi-local Noetherian ring; see~\cite[\S 34]{Ma} and~\cite[\S 7.8]{EGA IV}).

We introduce the following notations (see \cite{P2}).
Let $\delta = (p_0, \ldots, p_n)$.
Let a subscheme $X_i =  \overline{\{p_i\}}  $
be the closure of the point $p_i$ in $X$.
We introduce by induction the schemes $X_{i,\alpha_i}^{\prime}$
in the following  diagramm

$$
\begin{array}{ccccccc}
X_0           &  \supset  &  X_1                 & \supset &  X_2     & \supset  & \cdots  \\
\uparrow      &           &  \uparrow            &         &          &          &         \\
X'_0  &  \supset  &  X_{1, \alpha_1}     & \supset & \uparrow &          &         \\
              &           &  \uparrow            &         &          &          &         \\
              &           & X'_{1, \alpha_1} & \supset & X_{2, \alpha_2} &  &      \\
          &           &                          &         &  \uparrow       &  &      \\
          &           &                          &         &  \vdots         &  &      \\
          \end{array}
$$

Here $X'$ is the normalization of a scheme $x$
and $X_{i, \alpha_i}$
is an integral irreducible subscheme in $ X'_{i-1, \alpha_{i-1}}$
which is mapped onto $X_i$.
Then by any such diagram we obtain the collection of indices $ (\alpha_1, \ldots \alpha_n)$.
We denote the finite set of all such collections of indices by $\Lambda_{\delta}$.

Such a collection of indices $ (\alpha_1, \ldots \alpha_n) \in  \Lambda_{\delta}$
determines a chain of discrete valuations in the following way.
Integral irreducible subvariety $X_{1, \alpha_1}$ of the  normal scheme $X'_{0}$
defines the discrete valuation of the field of functions on $X_0$.
The residue field of this discrete valuation is the field of functions on the normal scheme
$X'_{1, \alpha_1}$, and  the integer irreducible subscheme $X_{2, \alpha_2}$
defiines the discrete valuation here. We have to proceed further for
$ \ldots, \alpha_3, \ldots, \alpha_n$ in this way.

Moreover, there is the following theorem (~\cite{P2}).(See also~\cite[theorem 3.3.2]{Y}
for the proof).

\begin{Th} \label{tl}
Let $X$
be an integral excellent $n$-dimensional Noetherian scheme.
Then for $\delta = (p_0, \ldots, p_n)$
the ring $K_{\delta}$ is an Artinian ring and
$$
K_{\delta} = \prod_{(\alpha_1, \ldots, \alpha_n) \in \Lambda_{\delta}}
K_{(\alpha_1, \ldots, \alpha_n)} \quad \mbox{,}
$$
where every $ K_{(\alpha_1, \ldots, \alpha_n)}  $
is an $n$-dimensional local field.
\end{Th}

\vspace{0.5cm}
\begin{ex} \em
To illustrate this theorem  we compute now
the ring $K_{\delta}$
in the following situation.
Let $p$ be a smooth point on an irreducible algebraic surface $X$ over $k$.
Suppose an irreducible curve $C \subset X$ contains the point $p$, but $C$
has at the point $p$ the node singularity, i.e., completed local ring  of the point $p$ on
the curve $C$ is $k[[t,u]]/tu$ for some local formal parameters
$u, t$ of the point $p$ on $X$. Let $\delta = (\eta, C, p)$,
where $\eta$ is a general point of $X$.
We fix any affine neighbourhood $\Spec B$ of $p$ on $X.$
Then according to formula~(\ref{f})
$$ C_{p} S_{p}^{-1} B = k(p) [[u,t]] $$
$$ C_{C}  S_{C}^{-1}  C_{p} S_{p}^{-1} B = k(p) ((u))[[t]] \oplus k(p) ((t))[[u]] $$
$$ K_{\delta} = C_{\eta}  S_{\eta}^{-1}
 C_{C}  S_{C}^{-1}  C_{p} S_{p}^{-1} B = k(p) ((u)) ((t)) \oplus  k(p) ((t))((u)) \mbox{.}
$$
\end{ex}

\section{Adeles and adelic complexes} \label{adel}
\subsection{Adeles on curves}
\label{adcur}

Let $C$ be a smooth connected algebraic curve over the field $k$.

For any  coherent sheaf $\f$ on $C$
we consider an adelic space $\da_C (\f)$:
$$
\da_C (\f) = \left\{ \{f_{p} \} \in \prod_{p \in C} \f \otimes_{\oo_C} K_{p}   \quad
\mbox{such that} \quad f_p \in    \f \otimes_{\oo_C} \oo_{K_{p}}
 \quad
 \mbox{for almost all}
  \quad
  p
\mbox{,}
\right\}
$$
where the product is over all closed points
 $p$
  of the curve   $C$.

We construct the following complex $ {\ad_C}(\f)$:
$$
\begin{array}{ccccc}
 \f \otimes_{\oo_C} k(C) & \times & \prod\limits_{p \in C}
 \f \otimes_{\oo_C} \oo_{K_{p}}  &
\lto  & \da_C (\f) \\
a & \times & b & \mapsto & a + b
 \mbox{.}
\end{array}
$$

We have the following theorem (for example,  see~\cite{S}).

\begin{Th}  \label{proposi}
The cohomology groups of the complex
 ${\ad_C}(\f)$
coincide with the cohomology $H^*(C, \f)$,
where $\f$ is
 any coherent sheaf  on  $C$.
\end{Th}

\proof We give here the sketch of the proof.
We write an adelic complex ${\ad_U}(\f)$ of the sheaf $\f$
for any open subset $U \subset C$.
Taking into account all $U$ we obtain a complex of sheaves $\ad (\f)$ on the curve $C$.
Then for the small affine $U$ we obtain that
the following complex
$$
0 \lto \f(U) \lto \ad_U (\f) \lto 0
$$
is exact, since we can apply the approximation theorem for Dedekind rings over fields.
Therefore the complex $\ad (\f)$ is a resolution of the sheaf $\f$ on $C$.
And by construction this resolution is a flasque resolution of the sheaf $\f$  on $C$.
Therefore it calculates the cohomology of the sheaf $\f$ on the curve $C$.

\subsection{Adeles on higherdimensional schemes} \label{adeles}

In this section we give a generalization of adelic complexes to the schemes
of arbitrary dimensions.

For algebraic surfaces adelic complexes were introduced by A.N. Parshin in~\cite{P1}.
We will give a detail exposition of adelic complexes on algebraic surfaces later
as an application of general machinery,  which was constructed
for arbitrary Noetherian schemes
by A.Beilinson in~\cite{B}. For a good exposition and proofs of Beilinson results
see~\cite{H}.

\subsubsection{Definition of adelic spaces.}
We introduce the following notations.
For any Noetherian scheme $X$ let $P(X)$ be the set of points of the scheme $X$.
Consider $p, q \in P(X)$. Define $p \ge q$
if $q \in \overline{\{ p\}}$, i.e., the point $p$ is in closure of the point $q$.
Then $\ge$ is a half ordering on $P(X)$.
Let $S(X)$ be the simplicial set induced by $(P(X), \ge)$,
i.e.,
$$S(X)_m = \{ (p_0, \ldots, p_m) \mid \nu_i \in P(X); p_i \ge p_{i+1} \}$$
is the set of $m$-simplices of $S(X)$
with the usual boundary  $\delta_i^n$ and degeneracy maps $\sigma_i^n$
for $n \in \dn$, $0 \le i \le n$.

Let $K \subset S(X)_n$.
For $p \in P(V)$ we denote
$$
{}_{p} K = \{ (p_1 > \ldots > p_{n}) \in S(V)_{n-1}
 \mid (p > p_1 \ldots  > p_n) \in K \} \mbox{.}
$$

Let ${\rm \bf QS}(X)$ and $ {\rm \bf CS}(X)$
be the category of  quasicoherent and coherent sheaves on the scheme $X$.
Let ${\rm  \bf Ab}$  be the category of Abelian groups.

We have the following proposition, see \cite{B}, \cite{H}, \cite{H1},
which is also a definition.

\begin{prop} \label{ind}
Let $S(X)$
be the simplicial set associated to the Noetherian scheme $X$.
Then there exist for integer $n > 0$, $K \in  S(X)_n$
functors
$$
\da(K, \cdot) \; : \;  {\rm \bf QS}(X)  \lto {\rm \bf Ab}
$$
uniquely determined by the properties \ref{pr1}, \ref{pr2}, \ref{pr3},
which are additive and exact.

\begin{enumerate}
\item  \label{pr1}
 $\da (K, \cdot)$
commutes with direct limits.
\item   \label{pr2}
For $n = 0$, a coherent sheaf $\f$ on $X$
$$
\da (K, \f) = \prod_{p \in K}
\mathop{\mathop{\lim}_{\leftarrow}}_l
\f_p / m_p^l \f_p  \mbox{.}
$$
\item \label{pr3}
For $n > 0$, a coherent sheaf $\f$ on $X$
$$
\da (K, \f) = \prod_{p \in P(V)}
\mathop{\mathop{\lim}_{\leftarrow}}_l
\da ( {}_{\eta} K, \f_p / m_p^l \f_p )
$$
\end{enumerate}
\end{prop}

\begin{nt} \em
Since any quasicoherent sheaf on an exellent Noetherian scheme
is a direct limit of coherent sheaves, we can apply  property~\ref{pr1}
of this proposition to define $\da (K, \f)$ on quasicoherent sheaves.
\end{nt}

\subsubsection{Local factors.}
The definition of $ \da (K, \f)$ is a definition of inductive kind.
By induction of the definition there is the following proposition, \cite[prop. 2.1.4.]{H}.
\begin{prop}
For integer $n > 0$, $K  \subset S(X)_n$, a  quasicoherent sheaf $\f$ on $X$
$$
\da (K, \f) \subset \prod_{\delta \in K} \da (\delta ,\f) \mbox{.}
$$
The inclusion is a transformation of functors.
\end{prop}

From this proposition we see that $\da (K, \f)$ is a kind of complicated adelic product inside of
$\prod\limits_{\delta \in K} \da (\delta,\f)$.

Therefore it is important to study the local factors $ \da (\delta,\f)$ for $\delta \in S(X)_n$. We have the following two propositions from \cite{H} about these local factors..

\begin{prop}
Let $\delta = (p_0, \ldots , p_n) \in S(X)_n$.
Let $U$ be an open affine subscheme which contains the point $p_n$
and therefore all of $\delta$. Let $M = \f(U)$. Then for a quasicoherent sheaf $\f$
$$
\da(\delta, \f) = \da (\delta, \tilde{M}) \mbox{,}
$$
where $\tilde{M}$ is a quasicoherent sheaf on  affine $U$ which corresponds to $M$.
\end{prop}

In the following proposition local factors $\da (\delta, \f)$
are computed for affine schemes.
\begin{prop}
Let $X = \Spec R$ and $\f = \tilde{M}$
for some $R$-module $M$. Further let $\delta \in (p_0, \ldots, p_n) \in S(X)_n$.
Then
\begin{equation} \label{for}
\da (\delta, \f) = C_{p_0} S_{p_0}^{-1} \ldots C_{p_n} S_{p_n}^{-1} R \otimes_R M
\mbox{.}
\end{equation}
$C_{p_0} S_{p_0}^{-1} \ldots C_{p_n} S_{p_n}^{-1} R$
is a flat Noetherian $R$-algebra. And for finitely generated $R$-modules
$$
C_{p_0} S_{p_0}^{-1} \ldots C_{p_n} S_{p_n}^{-1} R \otimes_R M =
C_{p_0} S_{p_0}^{-1} \ldots C_{p_n} S_{p_n}^{-1} M \mbox{.}
$$
\end{prop}

We compare now formula (\ref{for}) from the last proposition and
formula (\ref{f}) for $K_{\delta}$. We obtain that for an $n$-dimensional
Noetherian scheme $X$, for $\delta \in S(X)_n$ and a quasicoherent sheaf $\f$
$$
\da(\delta, \f) = K_{\delta} \otimes_{\oo_X} \f \mbox{.}
$$

\begin{nt} \em
 Due to  theorem \ref{tl} it means that for $\delta \in S(X)_n$ local factors $\da(\delta, \oo_X)$
on exellent Noetherian integral $n$-dimensional scheme $X$
are finite products of $n$-dimensional local fields.
\end{nt}

\subsubsection{Adelic complexes.}
\label{adcom}
Now we want to define  adelic complexes on the scheme $X$.

We have the simplicial set $S(X)$ with the usual boundary maps $\delta_i^n$
and degeneracy maps $\sigma_i^n$ for $n \in \dn$, $0 \le i \le n$.

We remark the following property, see~\cite[prop. 2.1.5.]{H}.

\begin{prop} \label{lp}
Let $K, L, M  \subset S(X)_n $
such that $K \cup M = L$, $K \cap M = \emptyset$.
Then there are natural transformations $i$ and $\pi$ of functors
$$ i(\cdot) : \da (K, \cdot) \lto \da (L, \cdot)
$$
$$
\pi(\cdot) : \da (L, \cdot)  \lto \da(M, \cdot)
$$
such that the following diagram is commutative and has split-exact lines for all
quasicoherent sheaves $\f$ on $X$
$$
\begin{array}{ccccccccc}
0 & \lto & \da (K, \f)  & \stackrel{i(\f)}{\lto} & \da (L, \f) &  \stackrel{\pi(\f)}{\lto}
& \da (M, \f) & \lto & 0
\\
& & \downarrow  & & \downarrow  & & \downarrow & \\
0 & \lto &  \prod\limits_{\delta \in K} \da (\delta, \f)  & {\lto} &
\prod\limits_{\delta \in L} \da (\delta, \f) &  {\lto}
& \prod\limits_{\delta \in M}  \da (\delta, \f) & \lto & 0 \mbox{.}
\end{array}
$$
\end{prop}

This proposition is proved by induction on definition-proposition~\ref{ind}.

\begin{defin}
Let $K \subset S(X)_0$, $\f$  a quasicoherent sheaf. Then let
$$
d^0 (K, \f)  \; : \;  \Gamma(X, \f)  \lto \da(K, \f)
$$
be the canoical map, which is a natural transformation of functors.
\end{defin}

\begin{defin} \label{ddef}
Let $K \subset S(X)_{n+1}$, $L \subset S(X)_n$,
$\delta_i^{n+1} K \subset L$ for some $i \in \{0, \ldots, n+1  \}$.
We define transformations of functors
$$
d_i^{n+1} (K, L, \cdot) \; : \; \da(L, \cdot)  \lto \da(K, \cdot)
$$
by the following properties.
\begin{enumerate}
\item If $i= 0$ and $\f$ is a coherent sheaf on $X$, then we apply the functor
$ \da( {}_{p} K, \cdot)$ to $\f \to \f_p/m_p^l \f_p $ and compose this map
with the projection of proposition~\ref{lp} for $L \supset {}_{p} K $.
We use the universal property of $\prod\limits_{p \in P(X)} \lim\limits_{\leftarrow}$.
\item
If $i = 1$, $n = 0$ and $\f$ is a coherent sheaf on $V$,
then the projection of proposition~\ref{lp} for $L \supset  {}_{p} K$
is composed with the following map. The maps $ d^0 ( {}_{p} K, \f_p / m_p^l \f_p)$
form a projective system for $l \in \dn$ and we apply
$\prod\limits_{p \in P(X)} \lim\limits_{\leftarrow}$ to it.
\item
If $i > 0$, $n >0$, $\f$ is a coherent sheaf, then the hypothesis $\delta_i^{n+1} K \subset L$
implies $\delta_{i-1}^n ({}_{p} K) \subset {}_{p} L$ for all $p \in P(X)$. Set
$$
d_i^{n+1} (K, L, \f) = \prod_{p \in P(X)}
\mathop{\lim\limits_{\leftarrow}}\limits_{l \in \sdbn}
d_{i-1}^n( {}_{p} K, {}_{p} L, \f_p / m_p^l \f_p)
\mbox{.}
$$
\item
$
d_i^{n+1} (K, L, \cdot)
$
commutes with direct limits.
\end{enumerate}
\end{defin}

\vspace{0.5cm}

For $\delta \in S(V)_{n+1}$
and $\delta' = \delta_i^{n+1} (\delta) \in S(V)_n $
by definition~\ref{ddef} we have local
boundary map
$$
d_i^{n+1}    \quad : \quad
\da(\delta' , \f)  \lto  \da(\delta , \f) \mbox{.}
$$
For $K \subset S(X)_{n+1}$, $L \subset S(X)_n$ with
$\delta_i^{n+1} K \subset L$
we define
$$
D_i^{n+1} (\f) \quad : \quad
\prod_{\delta \in L} \da(\delta, \f) \lto  \prod_{\delta \in K}   \da(\delta, \f ) \mbox{,}
$$
where $(x_{\delta})_{\delta \in L} \mapsto (y_{\delta})_{\delta \in K}$
is given by
$y_{\delta} = d_i^{n+1} (x_{\delta'})$.

For the computations of boundary maps  it is usefull the following proposition, which describes the boundary maps $d_i^{n+1}$ by means of the boundary maps $D_i^{n+1}$ on the product of local factors.

\begin{prop}
Let $K \subset S(X)_{n+1}$, $L \subset S(X)_n$
with $\delta_i^{n+1} K \subset L$. The following diagram commutes
$$
\begin{array}{ccc}
\da(L, \f)   & \stackrel{d_i^{n+1}}{\lto}  & \da (K, \f) \\
\downarrow   &                             & \downarrow \\
 \prod\limits_{\delta \in L} \da(\delta, \f) & \stackrel{D_i^{n+1}}{\lto}   &  \prod\limits_{\delta \in K}   \da(\delta, \f ) \mbox{.}
\end{array}
$$
\end{prop}

This porposition is proved by induction of definitions.

\vspace{0.5cm}

Now for the scheme $X$ we consider
the set $S(X)_n^{(red)}$ of non degenerate $n$-dimensional simplices.
(A simplex $(p_0, \ldots p_n)$ is  nondegenerate if $p_i \ne p_{i+1}$ for any $i$.)
For any $n \ge 0$, for any quasicoherent sheaf $\f$ on $X$ we denote
$$
\da^n_X (\f) = \da (S(X)_n^{(red)}, \f) \mbox{.}
$$
We consider the boundary maps
$$
d_i^{n+1} \quad  : \quad  \da^n_X  ( \f)  \lto \da^{n+1}_X (\f) \mbox{.}
$$

There are the following equalities for these boundary maps:
\begin{equation} \label{equ}
d_j^n d_i^n = d_i^n d_{j-1}^n   \qquad  \qquad i < j \mbox{.}
\end{equation}

For $n \ge 1$ we define $d_n : \da^{n-1}_X  (\f)  \lto \da^{n}_X (  \f)$ by
\begin{equation} \label{equ1}
  d^{n} = \sum_{j = 0}^{n} (-1)^j d_j^n \mbox{.}
\end{equation}

We  have the following proposition,
which is also a definition.
\begin{prop}
Differentials $d^n$ make $\da^*_X  ( \f)$
into a cohomological complex of Abelian groups $\ad_X (\f)$,
which we call the adelic complex of the sheaf $\f$ on $X$.
\end{prop}
\proof It follows by direct calculations with formulas~(\ref{equ}) and (\ref{equ1}).

We have the following theorem.
\begin{Th}
For any quasicoherent sheaf $\f$ on a Noetherian scheme $X$
$$
H^i (\ad_X (\f)) = H^i (X, \f) \mbox{.}
$$
\end{Th}
\proof
The proof of this theorem is a very far generalization
of the proof of theorem \ref{proposi}. Indeed, for any open subscheme $U \subset X$
we consider the following complex
\begin{equation} \label{cc}
0 \lto \f(U) \stackrel{d^0}{\lto} \da^0_U ( \f)  \stackrel{d^1}{\lto}  \da^1_U ( \f)
\stackrel{d^2}{\lto}  \ldots \stackrel{d^n}{\lto} \da^n_U ( \f)
\stackrel{d^{n+1}}{\lto} \ldots
\end{equation}
Taking into account all $U$, we obtain that this complex is a complex of sheaves
on $X$. Moreover,  by proposition~\ref{lp} the sheaves in this complex are flasque sheaves,
since $S(U)_n^{(red)}  \subset S(X)_n^{(red)}$ for any $n$.

By \cite[th.~4.1.1]{H}   for any affine scheme $U$ the complex~(\ref{cc})
is an exact complex. Therefore we constructed a flasque resolution of the sheaf $\f$ on $X$.
This resolution calculates the cohomology of the sheaf $\f$ on $X$.

\begin{nt} \em
We constructed here reduced adeles, since we used only nondegenerate simplices in $S(X)$.
These reduced adeles really carry information and they
are the part of the full complex, see \cite{H}.
\end{nt}

\subsection{Adeles on algebraic surfaces.}  \label{surface}
In this section we verify that the general adelic complex constructed in previous section
coincides with the adelic complex for curves from section~\ref{adcur}.
We give also application of general construction of adelic complexes
to algebraic surfaces.

We consider a smooth connected algebraic curve $C$ over  field $k$.
The set $S(C)_0^{(red)}$ consists of the general point $\eta$ and all closed points $p$ of the curve $C$. The set $S(C)_1^{(red)}$ consists of all pairs $(\eta, p)$.
For any coherent sheaf $\f$ on $C$ we can compute by definition
$$
\da_C^0 (\f) =
\f \otimes_{\oo_C} k(C) \; \times \; \prod_{p \in C}
 \f \otimes_{\oo_C} \oo_{K_{p}} \mbox{.}
$$

Let subset $K  \subset  S(C)_0^{(red)} $ consists of all closed points of the curve $C$.
We have by definition
$$
\da_C^1 (\f) = \da (K , \f_{\eta}) = \da (K , \f \otimes_{\oo_C} k(C)) =
$$
$$
= \da (K,  \mathop{\lim_{\lto}}_{D \in \Div(C)} \f \otimes_{\oo_C} \oo_C(D) =
$$
$$
= \mathop{\lim_{\lto}}_{D \in \Div(C)} \da (K, \f \otimes_{\oo_C} \oo_C(D) ) =
$$
$$
= \mathop{\lim_{\lto}}_{D \in \Div(C)}   \prod_{p \in C} \f \otimes_{\oo_C} \oo_{K_p}(D) \mbox{.}
$$

Therefore the adelic complex constructed in section~\ref{adcom}
coincides with the adelic complex for curves from section~\ref{adcur}.

\vspace{1cm}

We consider a smooth connected algebraic surface $X$ over field $k$.

The set $S(X)_0^{(red)}$
consists of the general point $\eta$ of $X$,
general points of all irreducible curves $C \subset X$,
all closed points $p \in X$.

The set $S(X)_1^{(red)}$ consists of all pairs $(\eta, C)$,
$(\eta, p)$ and $(C, p)$. (In our notations we identify
the general point of a curve $C \subset X$
with the curve $C$.)

The set $S(X)_2^{(red)}$ consists of all triples $(\eta, C, p)$.

\vspace{0.5cm}

We consider $\delta =  (\eta, C, p)$.
Let $f$ be a natural map from the local ring $\oo_{p,X}$  to the completion
$\hat{\oo}_{p,X}$. The curve $C$ defines a prime ideal $\C^{\,'}$ in the ring $\oo_{p,X}$.
Let $\C_1, \ldots, \C_n$ be all prime ideals of height $1$ in the ring $\hat{\oo}_{p,X}$
such that  for any $i$ $f^{-1} (\C_i) = \C^{\,'}$.
Any such $\C_i$ we will call a germ of $\C$ at $p$.
For any such germ $\C_i$ we define a two-dimensional
local field
$$
 K_{p,\C_i} =
  \Frac \;  {\mathop {\Lim_{\longleftarrow}}_l} \, (
 \mathop{\hoo_{p,X}}\nolimits_{ (\C_i ) }  / {\C}_i^{\, l}
\mathop{\hoo_{p,X}}\nolimits_{ (\C_i   ) } ) \mbox{.} $$

The ring $\mathop{\hoo_{x,X}}_{(\C_i)}$
is a localization of the ring
 $\hoo_{p,X}$
along the prime ideal $\C_i$.

Then according to formula~(\ref{f}),
we have (see~\cite{FP}, \cite{P1})
$$
\da (\delta, \oo_X) =  K_{\delta} = \bigoplus_{i=1}^{i=n} K_{p, \C_i} \mbox{.}
$$

Similarly  we have
$$
\da ((C, p), \oo_X ) = \bigoplus_{i=1}^{i=n} \oo_{ K_{p, \C_i}}  \mbox{.}
$$

We compute by definition
$$
\da ((\eta, C), \oo_X ) = K_C  \mbox{,}
$$
where the field $K_C$ is the completion of the field $k(X)$ along the discrete valuation
given by irreducible curve $C$ on $X$.

And from definition  we obtain that $
\da ((\eta, p), \oo_X)
$  is a subring in $\Frac(\hat{\oo}_{p, X}  )$
generated by subrings $k(X)$ and $\hat{\oo}_{p,X}$.
We denote this subring by $K_p$.

By definition we compute
$$
\da ((\eta), \oo_X) = k(X)  \mbox{,}
$$
$$
\da ((C), \oo_X) = \oo_{K_C}  \mbox{,}
$$
$$
\da ( (p), \oo_X ) = \hat{\oo}_{p,X} \mbox.
$$

\begin{nt}  \label{remar} \em
The local boundary  maps $d_i^n$ give natural embeddings of rings $\da ((\eta), \oo_X) $,
$\da ((C), \oo_X)$, $\da ( (p), \oo_X )$, $\da ((\eta, p), \oo_X)$,
 $\da ((\eta, C), \oo_X )$,
 $  \da ((C, p), \oo_X )$ to the ring $K_{\delta}$.
\end{nt}

By definition we have
$$
\da^0_X (\oo_X) \; = \;
k(X) \; \times \prod_{C \subset X} \oo_{K_C} \; \times \; \prod_{p \in X} \hat{\oo}_{p,X}
\mbox{.}
$$

\vspace{0.5cm}

From proposition~\ref{ind} and similarly to the case of algebraic curve we
compute the ring $ \da_X^2 (\oo_X) $, see details in~\cite{FP}, \cite{P2}.

For any prime ideal $\C \subset \hat{\oo}_{p, X}$ of height $1$ we define  the subring
$\hat{\oo}_{p,X} (\infty \C)$ of $ K_{p, \C}$:
$$\hat{\oo}_{p,X} (\infty \C) = \mathop{\lim\limits_{\lto}}\limits_l t_{\C}^{-l}
\hat{\oo}_{p,X} \mbox{,} $$
 where $t_{\C}$ is a generator
of ideal $\C$ in  $\hat{\oo}_{p,X}$.

The ring $\hat{\oo}_{p,X} (\infty \C)$ does not depend on
the choice of $t_{\C}$.

By $p \in \C \subset X$ we denote a germ at $p$ of an irreducible curve  on $X$.
Now we have
$$
\da_X^2 (\oo_X) \; =
\;  \{ f_{p,\C} \} \in \prod_{p \in \C \subset X} K_{p,\C} \qquad  \mbox{under the following  two conditions.}
$$
\begin{enumerate}
\item There exists a divisor $D$ on $X$ such that for any $p \in \C  \subset  X$
$$
\nu_{K_{p,\C}} ( f_{p,\C} ) \ge \nu_{\C} (D) \mbox{.}
$$
 \item For any irreducible curve $C \subset X$, any integer $k$ and all except
 a finite number  of points $p \in C$ we have that
inside of group $(K_{p,\C} \; \mod \C^{\, k} \oo_{K_{p, \C}})$
 $$
 f_{p,C} \; \mod \C^{\, k} \oo_{K_{p, \C}}  \; \in \;
 \hat{\oo}_{p,X} (\infty \C) \; \mod \C^{\, k} \oo_{K_{p, \C}} \mbox{.}
 $$
 Here we supposed that the curve $C$ has at $p$ a germ $\C$.
\end{enumerate}

\vspace{0.5cm}

We have
$$
\da^1_X (\oo_X) =  (\prod_{C \subset X} K_C) \cap  \da^2_X (\oo_X)  \; \times
 \; (\prod_{p \in X} K_p) \cap  \da^2_X (\oo_X) \;
\times \;
(\prod_{p \in \C \subset X }
\oo_{K_{p,\C}}) \cap  \da^2_X (\oo_X) \mbox{,}
$$
where we take  the intersection  inside of $\prod\limits_{p \in \C \subset X} K_{p, \C}$
due to remark \ref{remar} and diagonal embeddings
$$\prod\limits_{C \subset X} K_C  \lto
\prod\limits_{p \in \C \subset X } K_{p,\C} $$ and
$$\prod\limits_{p \in X} K_p \lto \prod\limits_{p \in \C \subset X } K_{p,\C} \mbox{.}$$

From formula \ref{equ1}  and explicit description of rings $\da^*_X(\oo_X)$ it is easy to
see differentials $d^n$ in the complex $\ad_X (\oo_X)$  ( \cite{FP}, \cite{P2}).
Indeed, let
$$
A_0 = k(X) \mbox{,} \qquad  A_1 = \prod_{C \subset X} \oo_{K_C} \mbox{,} \qquad
 A_2 =  \prod_{p \in X} \hat{\oo}_{p,X} \mbox{,}
$$
$$
A_{01} = (\prod_{C \subset X} K_C) \cap  \da^2_X (\oo_X) \mbox{,} \qquad \qquad
A_{02} =
(\prod_{p \in X} K_p) \cap  \da^2_X (\oo_X) \mbox{,}
$$
$$
A_{12} =
 (\prod_{p \in \C \subset X }
\oo_{K_{p,\C}}) \cap  \da^2_X (\oo_X)  \mbox{,}  \qquad \qquad
A_{012} = \da^2_X (\oo_X)
 \mbox{.}
$$
Then adelic complex $\ad_X (\oo_X)$ is
$$
\begin{array}{ccccc}
A_0 \oplus A_1 \oplus A_2 & \lto    & A_{01} \oplus A_{02} \oplus A_{12} &
\lto & A_{012} \\
(a_0, a_1, a_2)            & \mapsto & (a_1 - a_0, a_2 - a_0, a_2 -a_1)  &
& \\
                           &        & (a_{01}, a_{02}, a_{12}             &
               \mapsto & a_{01} - a_{02} + a_{12}  \mbox{.}
\end{array}
$$

\begin{nt} \label{interes}  \em
We remark the following interesting property,
see \cite[remark 5]{P5}, \cite{FP}.
 For any subset $I  \subset [0,1,2]$ we have
an embedding $A_I \hookrightarrow A_{012}$. Now for  any  subsets $ I, J  \subset [0,1,2]$
 we have that inside of group $A_{012}$
$$ A_I \cap A_J = A_{I \cap J} \mbox{.}
$$

This property is also true for corresponding components of adelic complex of any locally free sheaf on $X$.
\end{nt}

\section{Restricted adelic complexes} \label{sect4}
In this section we describe restricted adelic complexes.
The main difference of restricted adelic complexes from  adelic complexes constructed
in section~\ref{adel} is that restricted
adelic complexes are connected with one fixed chain (or flag) of irreducible subvarieties
of a scheme $X$.

Restricted adelic complexes come from so-called Krichever correspondence,~\cite{P5},~\cite{Os},
but see also~\cite{P4}
for connections with the theory of $\zeta$-functions of algebraic
curves.
Restricted adelic complexes on algebraic curves come originally from the
theory of integrable systems, see~\cite{SW}. For  algebraic surfaces restricted adelic complexes
were constructed by A.N. Parshin in \cite{P5}.
For higher dimensional schemes restricted adelic complexes were constructed by author
in~\cite{Os}.

\subsection{Restricted adelic complexes on algebraic curves and surfaces.} \label{curve}

We consider an irreducible algebraic curve $C$ over $k$.
We fix a smooth closed point $p \in C$.
For any coherent sheaf $\f$ of rank $r$ on $C$ we consider the following complex
\begin{equation} \label{cad}
\begin{array}{ccc}
\Gamma (C \setminus p, \f)  \;  \oplus \;
(\f \otimes_{\oo_C} \oo_{K_p} ) &
\lto & \f \otimes_{\oo_C} K_p  \\
(a_0 \oplus a_1) & \mapsto  &  a_1 -a_0
\mbox{.}
\end{array}
\end{equation}

We note that for the torsion free sheaf $\f$ we have natural embeddings
$$
\f \otimes_{\oo_C} \oo_{K_p}   \lto \f \otimes_{\oo_C} K_p
$$
$$
\Gamma (C \setminus p, \f) \lto \f \otimes_{\oo_C} K_p  \mbox{,}
$$
where the last embedding is given by
$$
\Gamma (C \setminus p, \f)  \lto
\Gamma (\Spec \oo_p \setminus p, \f) \lto
\Gamma (\Spec \oo_{K_p} \setminus p, \f) = \f {\otimes}_{\oo_C} K_p \mbox{.}
$$
Besides,  after the choice of basis of module $\f_p$ over the ring $\oo_p$ we have
$$
\f \otimes_{\oo_C} K_p = K_p^{\oplus r} \mbox{.}
$$
Therefore in this case complex~(\ref{cad}) is a complex of subgroups inside of $K_p^{\oplus r}$, where $K_p$ is a $1$-dimensional local field.

There is the following theorem, see, for example,  \cite{P5}, \cite{P4}.

\begin{Th} \label{th33}
The cohomology groups of complex~(\ref{cad})
coincide with the cohomology groups $H^*(C, \f)$.
\end{Th}

Chain of quasi-isomorphisms between complex~(\ref{cad}) and adelic complex $\ad_C(\f)$
was constructed   in~\cite{P5}. It proves the theorem~\ref{th33}.
We remark that it is important for the proof, that $C  \setminus p$
is an affine curve, see also remark~\ref{Chech} below.

The complex~(\ref{cad}) is called {\em restricted} adelic complex on $C$
associated with the point $p$.

\vspace{1cm}

Now let $X$  be an algebraic surface over $k$. We fix
 an irreducible  curve $C \subset X$,
and  a point $p \in C$ which is a smooth point on both $C$ and $X$.
Let $\f$ be a torsion free coherent sheaf on $X$.
We introduce the following notations from~\cite{P5}, \cite{P4}. Let $x \in C$,
$$
\hat{\f}_x  \qquad \mbox{,} \qquad \hat{\f}_C   \qquad \mbox{,} \qquad \hat{\f}_{\eta}
 $$
be completions of stalks of the sheaf $\f$ at scheme points given by  $x$ ,  irreducible curve $C$
and general point $\eta$ of $X$ correspondingly.
$$
B_x (\f) = \bigcap_{\D \ne C} (( \hat{\f}_x  \otimes K_x) \cap ( \hat{\f}_x \otimes \oo_{K_{x,\D}}) \mbox{,}
$$
where $\D$ runs over all germs at $x$ of irreducible curves
on $X$, which are not equal to $C$, and the intersection is done inside of the group $ \hat{\f}_x \otimes K_x$,
$$
B_C (\f) = (\hat{\f}_C \otimes K_C)  \cap \left(  \bigcap_{x \ne p} B_x  \right) \mbox{,}
$$
where the intersection is done inside of $\hat{\f}_p \otimes K_{x, \C}$
for all closed points $ x \ne p $ of $C$ and all germs $\C$  at $x$ of $C$,
$$
A_C (\f) = B_C (\f)  \cap \hat{\f}_C \mbox{,}
$$
$$
A(\f) =
 \hat{\f}_{\eta} \cap \left(  \bigcap_{x \in X -C}      \hat{\f}_p  \right)  \mbox{.}
$$

We note that
$$
A(\f) = \Gamma (X - C, \f)
$$
and for the smooth point $x \in C$  the space $B_x (\oo_X)$ coincides with the space
   $\hat{\oo}_{p,X} (\infty \C)$
   from section~\ref{surface}.

The following theorem was proved in \cite[th. 3]{P5}.
\begin{Th} \label{teres}
Let $X$ be an irreducible algebraic surface over a field $k$,
$C \subset X$ be an irreducible   curve,
and $p \in C$ be a smooth point on both $C$ and $X$.
Let $\f$ be a torsion free coherent sheaf on $X$.
Assume that the surface $X - C$ is affine. Then there exists a chain of quasi-isomorphisms
between  adelic complex $\ad_X ( \f)$ and the following
 complex
\begin{equation} \label{rescom}
\begin{array}{ccccc}
 A(\f) \oplus A_C (\f) \oplus  \hat{\f}_p & \lto    &  B_C (\f) \oplus  B_p (\f) \oplus
  (\hat{\f}_p \otimes {\oo}_{K_{p, \C}}) &
\lto &   \hat{\f}_p \otimes K_{p, \C} \\
(a_0, a_1, a_2)            & \mapsto & (a_1 - a_0, a_2 - a_0, a_2 -a_1)  &
& \\
                           &        & (a_{01}, a_{02}, a_{12}             &
               \mapsto & a_{01} - a_{02} + a_{12}  \mbox{.}
\end{array}
\end{equation}
\end{Th}

Under the conditions of this theorem the cohomology groups of the complex~(\ref{rescom})
coincide with the cohomology groups of adelic complex complex $\ad_X ( \f)$,
and therefore they are equal to $H^*(X, \f)$.

\begin{defin}
Complex~(\ref{rescom}) is called restricted adelic complex on $X$ associated with
the curve $C$ and the point $p \in C$.
\end{defin}

There is the following proposition, see~\cite[prop.4]{P5}.
\begin{prop} \label{ppp}
Under the conditions of theorem~\ref{teres}
we suppose also that $\f$ is a locally free sheaf, $X$ is a projective variety,  the local rings of $X$ are Cohen-Macaulay
and the curve $C$ is locally complete intersection. Then, inside the field
$K_{x, \C}$, we have
  $$
  B_C (\f) \cap B_p (\f) = A (\f)  \mbox{.}
  $$
\end{prop}

Let the rank of $\f$ be $r$.
Then after the choice of basis of $1$-dimensional free $\hat{\oo}_{p,X}$-module $\hat{\f_p}$
we have
$$
\hat{\f_p} = \hat{\oo}_{p,X}^{\oplus r} \mbox{,}
$$
$$
\hat{\f}_p \otimes K_{p, \C} = K_{p, \C}^{\oplus r} \mbox{,}
$$
$$
\hat{\f}_p \otimes \hat{\oo}_{K_{p, \C}} = \hat{\oo}_{K_{p, \C}}^{\oplus r} \mbox{,}
$$
$$
B_p(\f) = B_p ^{\oplus r} = \hat{\oo}_{p,X} (\infty \C)^{\oplus r} \mbox{,}
$$
$$
A_C (\f) = A (\f) \cap \hat{\oo}_{K_{p, \C}}^{\oplus r} \mbox{,}
$$
where the last intersection is done inside of $K_{p, \C}^{\oplus r}$.

Now due to proposition~\ref{ppp} we obtain that complex~(\ref{rescom})
is a complex of subgroups of $K_{p, \C}^{\oplus r}$
and is uniquely determined by one subgroup $B_C (\f)$ of $K_{p, \C}^{\oplus r}$.
In fact, all the other components of complex~(\ref{rescom}) can be defined
by intersections of $B_C (\f)$ with subgroups of $K_{p, \C}^{\oplus r}$,
which do not really depend on the sheaf $\f$.

\subsection{Restricted adelic complexes on higherdimensional schemes.}

In this section we construct restricted adelic complexes for arbitrary schemes.
These complexes will generalize corresponding complexes from section~\ref{curve}.

\subsubsection{General definitions.}
Let $X$ be
a Noetherian separated scheme.
Consider a flag of closed subschemes
$$
X \supset Y_0 \supset Y_1 \supset \ldots \supset Y_n
$$
in  $X$.
Let $J_j$ be the ideal sheaf of  $Y_j$ in $X$, $0 \le j \le n$.
Let $i_j$ be the embedding
 $Y_j  \hookrightarrow X$.
Let $U_i$ be an open subscheme of  $Y_i$
complementing $Y_{i+1}$, $0 \le i \le n-1$.
Let $j_i : U_i \hookrightarrow Y_i$ be the open embedding
of  $U_i$
in $Y_i$, $0 \le i \le n-1$.
Put
 $U_n = Y_n$
and let $j_n$ be
the identity morphism from
 $U_n$  to $Y_n$.

Assume that every  point
 $x \in X$
has an open affine neighbourhood
 $U \ni x$
such that
$U \cap U_i$ is an affine scheme for any $0 \le i \le n$.
In what follows, a flag  of subschemes
 $\{Y_i, \; 0 \le i \le n \}$ with this condition
is called a flag with { \em locally affine complements.}

\begin{nt} { \em
The last condition (existence of locally affine complements) holds, for example,
in the following cases:
\begin{itemize}
\item $Y_{i+1}$
is the Cartier divisor on
$Y_i$ for $0 \le i \le n-1$), and
\item
$U_i$ is an affine scheme for any  $0 \le i \le n-1$ (the intersection of two open affine subschemes on a separated scheme
is an  affine subscheme.)
\end{itemize}
}
\end{nt}

\bigskip
Consider the $n$-dimensional
simplex
and its standard simplicial set (without degeneracy).
To be precise, consider the set:
$$
(\{ 0\}, \{ 1\}, \ldots, \{ n\})
$$
(all the integers between $0$ and $n$.)
Then the simplicial set $S = \{ S_k \}$ is given by
\begin{itemize}
\item $S_0 \eqdef \{\eta \in \{ 0\}, \{ 1\}, \ldots, \{ n\}   \} $.
\item
$ S_k \eqdef \{ (\eta_0, \ldots, \eta_k),
\quad \mbox{where} \quad \eta_l \in S_0  \quad \mbox{and} \quad
\eta_{l-1} < \eta_l  \} $.
\end{itemize}
The boundary map $\partial_i$ ($0 < i < k$)
is given by eliminating the
$i$-th component of the vector  $(\eta_0, \ldots, \eta_k)$
to give
 the  $i$-th face  of  $(\eta_0, \ldots, \eta_k)$.

Let $ {\rm \bf QS}(X)$ be the category  of quasicoherent sheaves on $X$.
Let ${\rm \bf Sh}(X) $ be the category of sheaves of Abelian groups on $X$.
Let $f : Y \longrightarrow X$ be a morphism of  schemes.
Then $f^*$ always denotes
the pull-back   functor  in the category
of sheaves of Abelian groups, and $f_*$ is
the direct image functor
in the category of sheaves of Abelian groups.

We give the following definition from~\cite{Os}.
\begin{defin}
For any  $(\eta_0, \ldots , \eta_k)  \in S_k$
we define a functor
$$
V_{(\eta_0, \ldots, \eta_k)} \; : \;  {\rm \bf QS}(X) \lto {\rm \bf Sh}(X)
\mbox{,}
$$
which is  uniquely determined by the following inductive conditions:
\begin{enumerate}
\item
$V_{(\eta_0, \ldots, \eta_k)}  $
commutes with direct limits.
\item
If $\f$ is a coherent sheaf and $\eta \in S_0$,
then
$$
V_{\eta}(\f) \eqdef
\mathop{\pl}\limits_{m \in \bf{N}}
(i_{\eta})_* (j_{\eta})_*
(j_{\eta})^* (\f / J^m_{\eta} \f) \mbox{.}
$$
\item
If $\f$ is a coherent sheaf
and $(\eta_0, \ldots, \eta_k) \in S_k$, $k \ge 1$,
then
$$
V_{(\eta_0, \eta_1, \ldots, \eta_k)}(\f) \eqdef
\mathop{\pl}\limits_{m \in {\bf N}}
V_{(\eta_1, \ldots, \eta_k)}
\left( (i_{\eta_0})_*  (j_{\eta_0})_*  (j_{\eta_0})^*
(\f / J^m_{\eta_0} \f) \right)  \mbox{.}
$$
\end{enumerate}
\end{defin}

 We'll use sometimes
the equivalent  notation for
$V_{(\eta_0, \ldots, \eta_k)}(\f)$,
in which the closed subschemes are indicated explicitly:
$$
V_{(\eta_0, \ldots, \eta_k)}(\f) =
V_{(Y_{\eta_0}, \ldots, Y_{\eta_k})}(X, \f) \mbox{.}
$$

There is
the following proposition, \cite[prop. 1]{Os}, which is proved by induction.
\begin{prop} \label{predl1}
Let $ \sigma = (\eta_0, \ldots, \eta_k) \in S_k$.
Then the following assertions hold.
\begin{enumerate}
\item       \label{pun1}
The functor $V_{\sigma} :
\QS(X) \lto \Sh(X)$ is well defined.
\item        \label{pun2}
The functor $V_{\sigma}$
is exact and additive.
\item     \label{pun3}
The functor $V_{\sigma}$
is local on $X$, that is,
for any open
$U \subset X$
and any quasicoherent sheaf
 $\f$ on $X$ we have
$$
V_{(Y_{\eta_0}, \ldots, Y_{\eta_k})}(X, \f)  \mid_U =
V_{(Y_{\eta_0} \cap U, \ldots, Y_{\eta_k} \cap U  )} (U, \f \mid_U)  \mbox{.}
$$
(if  $Y_j \cap U = \o$,
then $Y_i \cap U$ is the empty subscheme of $U$
defined by the ideal sheaf $\oo_U$.)
\item For any quasicoherent sheaf $\f $ on $X$
the sheaf $V_{(\eta_0, \ldots, \eta_k)}(\f)$ is a sheaf of
 $\oo_X$-modules
 supported  on the subscheme $Y_{\eta_k}$.
(In general, this sheaf  is not quasicoherent.)
\item   \label{pun5}
For any quasicoherent sheaf
 $\f$ on $X$ we have
$$
V_{\sigma}(\f) =
V_{\sigma}(\oo_X)  \otimes_{\oo_X} \f \mbox{.}
$$
\item    \label{pun7}
If all $U_i$ is affine, $0 \le i \le n$,
then for any quasicoherent
sheaf $\f$ on $X$ and any $m \ge 1$
we have
$$
H^m (X, V_{\sigma}(\f)) = 0    \mbox{.}
$$
\end{enumerate}
\end{prop}

\vspace{0.5cm}

\subsubsection{Construction of restricted adelic complex.}
We consider  the standard $n$-simplex
$S= \{ S_k, \; 0 \le k \le n \}$ without degeneracy.
If $\sigma= (\eta_0, \ldots, \eta_k) \in S_k$,
then $\partial_i (\sigma)$ is the $i$th face of $\sigma$, $0 \le i \le k$.
We {\em define} a morphism of functors
$
d_i(\sigma) \quad : \quad  V_{\partial_i(\sigma)} \lto V_{\sigma}$,
as the morphism that
commutes with direct limits and
coincides
on coherent sheaves with the map
\begin{equation}   \label{tank}
V_{\partial_i(\sigma)}(\f) \lto V_{\sigma}(\f)  \mbox{}
\end{equation}
defined by the following rules.
\begin{itemize}
\item[a)]
If $i =0$,
then (\ref{tank})
is obtained by applying  the functor
$V_{\partial_0 (\sigma)}$
to the map
$$
\f \lto
(i_{\eta_0})_* (j_{\eta_0})_*  (j_{\eta_0})^*
(\f / J_{\eta_0}^m \f)
$$
and passing to the projective limit with respect to $m$;
\item[b)]
If $i=1$ and $k=1$,
then we have the natural map
$$
(i_{\eta_0})_* (j_{\eta_0})_*  (j_{\eta_0})^*
(\f / J_{\eta_0}^m \f)
\lto
V_{(\eta_1)} ((i_{\eta_0})_* (j_{\eta_0})_*  (j_{\eta_0})^*
(\f / J_{\eta_0}^m \f)) \mbox{.}
$$
Passing to the projective limit with respect to
 $m$,
we get the map~(\ref{tank}) in this case.
\item[c)]
If  $i \ne 0$ and $k >1$,
then  we use induction on
$k$ to get the map
$$
V_{\partial_{i-1} \cdot
(\partial_0 (\sigma))
} ((i_{\eta_0})_* (j_{\eta_0})_*  (j_{\eta_0})^*
(\f / J_{\eta_0}^m \f))
\lto
V_{
\partial_0 (\sigma)
} ((i_{\eta_0})_* (j_{\eta_0})_*  (j_{\eta_0})^*
(\f / J_{\eta_0}^m \f))  \mbox{.}
$$
Passing to the projective limit to $m$
we get the map~(\ref{tank})
in this case.
\end{itemize}

There is the following proposition,~\cite[prop. 3]{Os}.
\begin{prop}     \label{predl3}
For any $1 \le k \le n$, $0 \le i \le k$ let
$$
d_i^k  \eqdef \sum_{\sigma \in S_k} d_i(\sigma)
\quad : \quad \bigoplus_{\sigma \in S_{k-1}} V_{\sigma}
\lto
\bigoplus_{\sigma \in S_k} V_{\sigma}  \mbox{.}
$$
We define the map
$$d_0^0 \quad : \quad id \lto \bigoplus_{\sigma \in S_0}
V_{\sigma} \mbox{}$$
as the direct sum of the natural maps
$\f \lto
V_{\sigma} (\f)$. (Here $id$ is the functor
of the natural imbedding of
 $\QS(X)$ into $\Sh(X)$,
$\f$ is a quasicoherent sheaf on $X$,
$\sigma \in S_0$.)

Then for all $0 \le i < j \le k \le n-1$ we have
\begin{equation}  \label{npp}
 d_j^{k+1} d_i^k  = d_i^{k+1} d_{j-1}^k   \mbox{.}
\end{equation}
\end{prop}

Let
$$
d^m \eqdef \sum\limits_{0 \le i \le m} (-1)^i d^m_i
$$
Then,
given any quasicoherent sheaf $\f$ on $X$,
 proposition~\ref{predl3}  enables us
to construct
the complex of sheaves  $V(\f)$  in the  standard way:
$$
\ldots
\lto
\bigoplus_{\sigma \in S_{m-1}} V_{\sigma}(\f)
\stackrel{d^m}{\lto}
\bigoplus_{\sigma \in S_m} V_{\sigma}(\f)
\lto
\ldots \mbox{.}
$$
The property  $d^{m+1} d^m = 0$ follows
from~(\ref{npp})
by an easy  direct calculation.

\medskip
We have the following theorem,~\cite[th. 1]{Os}.
\begin{Th}   \label{teorem1}
Let $X$
be a Noetherian separated scheme and let
$Y_0 \supset Y_1 \supset \ldots  \supset  Y_n$ be a flag
of closed subschemes with locally affine complements.
Assumee that $Y_0 = X$.
Then the following complex is exact:
\begin{equation}  \label{kff}
0 \lto \f \stackrel{d^0}{\lto} V(\f) \lto 0     \mbox{.}
\end{equation}
\end{Th}
\proof. We give the sketch of the proof. It suffices to consider the case when the sheaf $\f$
is coherent. We consider the exact sequence of sheaves
$$
0 \lto \h \lto \f \lto (j_0)_* (j_0)^* \f \lto \g \lto 0 \mbox{.}
$$
Since the functors $V_{\sigma}$ are exact for all $\sigma$,
we obtain the following exact sequence of complexes of sheaves:
\begin{equation} \label{posledova}
0 \lto  V(\h) \lto V(\f) \lto V((j_0)_* (j_0)^* \f) \lto V(\g) \lto 0 \mbox{.}
\end{equation}
The sheaves $\h$ and $\g$ are supported on $Y_1$.
Therefore by induction we may assume that the complexes
$$
0 \lto \h \stackrel{d^0}{\lto} V (\h) \lto 0 \mbox{,}
$$
$$
0 \lto \g \stackrel{d^0}{\lto} V (\g) \lto 0
$$
are already exact.  The complex
$$ 0 \lto (j_0)_*(j_0)^*\f  \stackrel{d^0}{\lto}  V( (j_0)_* (j_0)^* \f ) \lto 0  $$
is exact,  because  the complex $V (j_0)_* (j_0)^* \f$ has the same components
$V_{\sigma'} (\f)$ of degrees $k$ and $k+1$, $\sigma' = (0, \eta_0, \ldots, \eta_k)$
for $\sigma = (\eta_0, \ldots, \eta_k) \in S_k$. Now the theorem follows from
sequence~(\ref{posledova}).

\vspace{0.5cm}
For any $\sigma \in S_k$ we define
$$
A_{\sigma}(\f) \eqdef  H^0(X, V_{\sigma} (\f))  \mbox{.}
$$

We have the following proposition,~\cite[prop. 4]{Os}.
\begin{prop} \label{predl4}
Let $X$ be a Noetherian separated scheme,
and let
$Y_0 \supset Y_1 \supset \ldots  \supset  Y_n$ be a flag
of closed subschemes
such that
$U_i$ is affine, $0 \le i \le n$.
Let $\sigma \in S_k$ be arbitrary.
\begin{enumerate}
\item        \label{punkt1}
Then
$A_{\sigma}$ is an exact and additive
functor: $\QS(X) \lto \Ab$.
\item  If $X = \Spec A$ and
$M$ is some $A$-module,
then
$$
A_{\sigma}(\tilde{M}) = A_{\sigma} (\oo_X)  \otimes_A M  \mbox{.}
$$
\end{enumerate}
\end{prop}

Let $\f$ be any quasicoherent sheaf on  $X$.
Applying the functor
 $H^0(X, \cdot)$ to the complex $V(\f)$,
we obtain the complex  $A(\f)$ of Abelian groups:
$$
\ldots \lto \bigoplus_{\sigma \in S_{m-1}} A_{\sigma}(\f)
\lto
\bigoplus_{\sigma \in S_{m}} A_{\sigma}(\f)
\lto
\ldots   \mbox{.}
$$

\medskip
Now we have the following theorem, see~\cite[th. 2]{Os}.
\begin{Th}    \label{teorem2}
Let $X$ be a Noetherian separated scheme.  Let
$Y_0 \supset Y_1 \supset \ldots \supset Y_n$ be a flag
of closed subschemes
such that
 $Y_0 = X$ and
$U_i$ is affine, $0 \le i \le n$.
Then the cohomology of the  complex
 $A(\f)$
coincide with that of the  sheaf
$\f$  on $X$,
that is,  for any $i$
$$
H^i(X, \f) = H^i(A(\f))   \mbox{.}
$$
\end{Th}
\proof
It follows from theorem~\ref{teorem1}
and assertion~\ref{pun7}  of proposition~\ref{predl1}
 that $V(\f)$ is
an acyclic resolution for the sheaf
$\f$.
Hence
the
cohomology   of
 $\f$
 may be calculated
by means of global sections of this resolution.
Theorem~\ref{teorem2} is proved. \\
\medskip \\
This theorem
immediatly yields the following geometric corollary, see~\cite[th. 3]{Os}.
\begin{Th}     \label{teorem3}
Let $X$ be a projective algebraic scheme
of dimension $n$ over a field.
Let
$Y_0 \supset Y_1 \supset \ldots \supset Y_n$ be
a flag of closed subschemes
such that
 $Y_0 = X$ and
$Y_i$
is an ample divisor on
$Y_{i-1}$ for
$1 \le i \le n$.
Then for any quasicoherent sheaf
 $\f$
on $X$ and
  any $i$  we have
$$
H^i(X, \f) = H^i(A(\f))   \mbox{.}
$$
\end{Th}
\proof Since
$Y_i$ is an ample divisor on $Y_{i-1}$
for  $1 \le i \le n$, we see that
 $U_i$ is an affine scheme for all $0 \le i \le n-1$.
Since $ \dm  Y_n =0$, we obtain that $U_n= Y_n$ is also affine.
Applying theorem~\ref{teorem2}
we complete the proof.

\smallskip
\medskip
\begin{nt} \label{Chech}
{\em
For any quasicoherent sheaf
 $\f$
and any $\sigma = (\eta_0) \in S_0$,
$A_{\sigma} (\f)$
is the group of section over $U_{\eta_0}$
of the sheaf $\f$
lifted to the formal neighbourhood
of the subscheme $Y_{\eta_0}$
in $X$.
The complex
 $A(\f)$ can be  interpreted as
the \v{C}ech complex for this ''acyclic covering ''
of the scheme $X$.
}
\end{nt}

\begin{defin}
The complex $A (\f)$
is called  restricted adelic complex on $X$ associated with flag
$
Y_0 \supset Y_1 \supset \ldots \supset Y_n
$.
\end{defin}

\begin{nt} { \em
\begin{itemize}
\item If $C$ is an algebraic curve, $Y_0 = C$ and $Y_1 = p$ is a  smooth point ,
then $A(\f)$ coincides with  complex~(\ref{cad}). Indeed,
$$
A_0(\f) = \Gamma (C \setminus p, \f) \mbox{,}
$$
$$
A_1 (\f) = \f \otimes_{\oo_C} \oo_{K_p}  \mbox{,}
$$
$$
A_{12} (\f) = \f \otimes_{\oo_C} K_p    \mbox{.}
$$
\item If $X$ is an algebraic surface, $Y_0 = C$ and $Y_1 = p$ is a smooth point
on both $C$ and $X$, then $ A(\f)$ coincides with  complex~(\ref{rescom}). Indeed,
$$
A_0 (\f) =    A(\f)   \mbox{,}
\qquad \qquad
A_1 (\f) = A_C (\f) \mbox{,}
$$
$$
A_2 (\f) = \hat{\f}_p  \mbox{,}
\qquad \qquad
A_{01} (\f) = B_C (\f) \mbox{,}
$$
$$
A_{02}  (\f) = B_p (\f) \mbox{,}
\qquad \qquad
A_{12}   (\f) = \hat{\f}_p \otimes {\oo}_{K_{p, \C}} \mbox{,}
$$
$$
A_{012} (\f) = \hat{\f}_p \otimes K_{p, \C} \mbox{.}
$$
\end{itemize}
}
\end{nt}

\medskip

\subsubsection{Reconstruction of restricted adelic complex.}
By propostition~\ref{predl3} we have
the natural map
$$
d_i (\sigma)  \: : \:
A_{\partial_i (\sigma)} (\f)  \lto A_{\sigma} (\f)
$$
for any $\sigma \in S_k$, $1 \le k \le n$,
and any $i$, $0 \le i \le k$.

Taking~(\ref{npp}) into account,
we obtain
the natural map
$$
\quad  A_{\sigma_1}(\f)  \lto A_{\sigma_2}(\f)  \quad
$$
for any locally free sheaf $\f$
on $X$ and any $\sigma_1, \sigma_2 \in S$, $ \sigma_1 \subset \sigma_2$.

There is the following proposition, see~\cite[th. 4]{Os}.
\begin{prop}      \label{tend1}
Let $X$
be a projective equidimensional
Cohen-Macaluay scheme of dimension
$n$  over a field.
Let $Y_0 \supset Y_1 \supset \ldots
\supset Y_n$ be a flag of closed subschemes
such that
 $Y_0 = X$ and
$Y_i$         is an ample Cartier divisor on
$Y_{i-1}$
for
 $1 \le i \le n$.
Then
the following assertions hold for any locally free sheaf $\f$ on $X$.
\begin{enumerate}
\item     \label{unkt1}
The natural map
$
H^0(X, \f) \lto A_{\sigma} (\f)
$
is an embedding
for any $\sigma  \in S_k$, $0 \le k \le n$.
\item    \label{unkt2}
The natural map
$
\quad  A_{\sigma_1}(\f)  \lto A_{\sigma_2}(\f)  \quad
$
is an embedding
for any locally free sheaf $\f$
on $X$ and any $\sigma_1, \sigma_2 \in S$, $ \sigma_1 \subset \sigma_2$.
\end{enumerate}
\end{prop}

\begin{nt} { \em
We note that any integral Noetherian scheme of dimension $1$
is a Cohen-Macaulay scheme.
Any normal Noetherian scheme of dimension $2$ is a Cohen-Macaulay scheme,
see~\cite[Ch. II, th. 8.22A]{Ha}.
}
\end{nt}

\bigskip

By $(0,1, \ldots, n) \in S$ we denote the unique face of dimension $n$ in $S$.

By proposition~\ref{tend1}
we can  embed
$A_{\sigma_1} (\f)$
and $A_{\sigma_2} (\f)$
in $A_{(0,1 \ldots, n)} (\f)$
for any $\sigma_1, \sigma_2 \in S$
and any locally free sheaf $\f$.

Now we formulate the following theorem, see~\cite[th. 5]{Os}.
\begin{Th}       \label{tend2}
Let all the hypotheses of proposition~\ref{tend1}
be satisfied.
Then the following assertions hold for any locally free sheaf $\f$
and any
$\sigma_1, \sigma_2 \in S$.
\begin{enumerate}
\item                \label{pu1}
If $\sigma_1 \cap \sigma_2 = \o$,
then, inside of $A_{(0,1, \ldots, n)}$,
$$
A_{\sigma_1} (\f)  \cap A_{\sigma_2} (\f) = H^0 (X, \f)  \mbox{.}
$$
\item                 \label{pu2}
If $\sigma_1 \cap \sigma_2 \ne \o$,
then, inside of $A_{(0,1, \ldots, n)}$,
$$
A_{\sigma_1}(\f) \cap A_{\sigma_2} (\f)
= A_{\sigma_1 \cap \sigma_2} (\f)  \mbox{.}
$$
\end{enumerate}

\end{Th}

\begin{nt} {\em
Theorem~\ref{tend2} is similar to the property of adelic complexes $\ad_X (\f)$ noticed
in remark~\ref{interes}.
}
\end{nt}

\vspace{0.5cm}

We  assume
that the hypotheses of proposition~\ref{tend1}
hold and that
the field $k$  of definition of the scheme $X$ is $k$.
We also  assume that
$Y_n=p$, where $p$ is a smooth point on each $Y_i$, $0 \le i \le n$.

Let us choose and fix  local parameters
 $t_1, \ldots, t_n \in
\widehat{\oo}_{p,X}$
such that
$t_{i} {|}_{Y_{i-1}} = 0$
is a local equation of the divisor
 $Y_i$
in the formal neighbourhood of the point~$p$
on the scheme $Y_{i-1}$, $1 \le
i \le n$.

Let $\f$ be a rank~$1$ locally free
sheaf on $X$.
We fix a trivialization $e_p$ of  $\f$
in a  formal neghbourhood   of the point
$p$ on $X$, that is an isomorphism
$$
e_p \quad : \quad \hat{\f}_p \lto \hat{\oo}_{p, X} \mbox{.}
$$

By the choice of the local parameters and the trivialization
we can identify
 $A_{(0,1, \ldots, n)}(\f)$
with the  $n$-dimensional local field  $k(p)((t_n))\ldots (((t_1))$:
$$
A_{(0,1, \ldots, n)}(\f) = k(p)((t_n))\ldots (((t_1)) \mbox{.}
$$

Moreover,
we fix a set of integers $0 \le j_1 \le \ldots \le j_k \le n-1$.
Define $\sigma_{(j_1, \ldots, j_k)} \in S_{n-k}$ as the following set
$$ \left\{i \: :  \: 0 \le i \le n, \, i \ne j_1, \, \ldots,
\,
i \ne j_k \right\} \mbox{.} $$
By proposition~\ref{tend1}  we have a natural embedding
$$A_{\sigma_{(j_1, \ldots, j_k)}} (\f) \lto  A_{(0,1, \ldots, n)}(\f)  \mbox{.}$$
When we identify
$A_{(0,1, \ldots, n)}(\f)$ with the field
$k(p)((t_n)) \ldots ((t_1))$,
the space
 $A_{\sigma_{(j_1, \ldots, j_k)}} (\f)$
corresponds to the following
 $k$-subspace in
$k(p)((t_n)) \ldots ((t_1))$:
\begin{equation}  \label{st}
\left\{  \sum a_{i_1, \ldots, i_n}  t_n^{i_n}
\ldots  t_1^{i_1 }  \; :  \;   a_{i_1, \ldots, i_n}  \in k(p), \,
i_{j_1+1} \ge 0,  \, i_{j_2+1} \ge 0, \, \ldots, i_{j_k +1} \ge 0
 \right\} \mbox{.}
\end{equation}

Thus, by theorem~\ref{tend2},
to  determine  the images of
$A_{\sigma} (\f)$
in $k(p)((t_n)) \ldots ((t_1))$
(for any $\sigma \in S$),
it suffices to know  only one image of
$A_{(0,1, \ldots, n-1)}$
in $k(p)((t_n)) \ldots ((t_1))$.
(All the others are obtained  by taking
the intersection of the image of this one
with the standard subspaces~(\ref{st})
of $k(p)((t_n)) \ldots ((t_1))$.)

It is clear that these arguments
are generalized immediately to locally free sheaves
 $\f$  of rank $r$
and to the spaces $k(p)((t_n)) \ldots ((t_1))^{\oplus r}$.

These arguments lead to the following theorem, which enables to reconstruct restricted adelic complex $A (\f)$, see also~\cite[th. 6]{Os}.
\begin{Th}
Let all the hypotheses of proposition~\ref{tend1}
be satisfied.
We also  assume that
$Y_n=p$, where $p$ is a smooth point on each $Y_i$, $0 \le i \le n$.
Let $\f$ be a locally free sheaf on $X$.
Then the subspace
$$
 A_{(0,1, \ldots, n-1)} (\f)
 \;
 \subset
 \;
 A_{(0,1, \ldots, n)}(\f)
  $$
 uniquely determines the restricted adelic complex $A (\f)$.
\end{Th}

\section{Reciprocity laws.}

Let $L$ be a  field of discrete valuation $\nu_L$ with the valuation ring $\oo_L$
and the maximal idela $m_L$.
Then there is the tame symbol
 \begin{equation} \label{tame}
(f,g)_L  = (-1)^{\nu_L(f)\nu_L(g)}
\frac{f^{\nu_L(g)}}{g^{\nu_L(f)}} \; \mod  \; m_L \mbox{,}
\end{equation}
 where $f$, $g$ are from $L^*$.

 There is the following reciprocity law, see, for example, \cite{S}.
\begin{prop}
Let $C$ be a complete smooth algebraic curve over $k$.
For any $f, g \in k(C)^*$ in the following product only a finite number
of terms is not equal to $1$ and
\begin{equation} \label{curver}
\prod_{p \in C} \nm\nolimits_{k(p)/k} \, (f, g)_{K_p} =1 \mbox{,}
\end{equation}
where the product is taken over all closed points $p \in C$,
and $\nm$ is the norm map.
\end{prop}

Let $K $ be a $2$-dimensional local field with the last residue field $k$.
 We have  a discrete valuation of rank~$2$ on $K$:
 $$(\nu_1, \nu_2) :  K^* \to \dz \oplus \dz  \mbox{.}$$
 Here $\nu_1 = \nu_K$ is the discrete
valuation of the filed $K$, and
$$\nu_2 (b) \eqdef \nu_{\bar{K}} ( \: \overline{b t_1^{-\nu_1(b)}} \: ) \mbox{,} $$
where $\nu_1 (t_1) =1$.
We note that $\nu_2$ depends on the choice of local parameter $t_1$.

Let $m_K$ be the maximal ideal of $\oo_K$, and
$m_{\bar{K}}$ be  the  maximal ideal of $\oo_{\bar{K}}$.

We define a map:
$$
\nu_K (\; , \;)  \quad : \quad  K^* \times K^* \lto \dz
$$
as the composition of maps:
$$
K^* \times K^*  \lto K_2(K) \stackrel{\partial_2}{\lto} \bar{K}^*  \stackrel{\partial_1}{\lto} \dz \mbox{,}
$$
where $\partial_i$ is the boundary map in algebraic $K$-theory. The map $\partial_2$ coincides with tame symbol~(\ref{tame}) with respect to
discrete valuation $\nu_1$. The map $\partial_1$ coincides with the discrete valuation $\nu_{\bar{K}}$.

We define a map:
$$
(\;, \;, \;)_K
\quad : \quad
K^* \times K^* \times K^* \lto k^*
$$
as the composition of maps
$$
K^* \times K^* \times K^* \lto  K_3^M (K) \stackrel{\partial_3}{\lto} K_2(\bar{K})
\stackrel{\partial_2}{\lto} k^* \mbox{,}
$$
where
$K_3^M$ is the Milnor $K$-group.

There are the following explicit expressions for these maps (see~\cite{FP}):
$$
\nu_K (f, g) = \nu_1(f) \nu_2(g) - \nu_2(f) \nu_1(g)
$$
$$
(f,g,h)_K = \sign\nolimits_K(f,g,h)
 f^{\nu_K (g, h)} g^{\nu_K (h, f)} h^{\nu_K (f, g)} \mod m_K
\mod\nolimits m_{\bar{K}}
$$
$$
\sign\nolimits_K(f,g,h) = (-1)^B  \mbox{,}
$$
where $
B =
\nu_1(f) \nu_2(g) \nu_2(h)
+ \nu_1(g) \nu_2(f) \nu_2(h)
+ \nu_1(h) \nu_2(g) \nu_2(f)
+ \nu_2(f) \nu_1(g) \nu_1(h)
+ \nu_2(g) \nu_1(f) \nu_1(h)
+ \nu_2(h) \nu_1(f) \nu_1(g)
$.

\begin{prop}
For any $f,g,h \in K^*$
$$
\sign\nolimits_K(f,g,h) =
(-1)^A \qquad \mbox{, where}
$$
$$
A = {\nu_K (f,g) \nu_K (f,h)
+ \nu_K (f,g) \nu_K (f,h)
+ \nu_K (f,g) \nu_K (f,h)
+ \nu_K (f,g) \nu_K (f,h) \nu_K (g,h)
} \mbox{.}
$$
\end{prop}
\proof It follows from direct calculations modulo $2$ with $A$ and $B$  using the explicit expressions above.

\vspace{0.7cm}

Let $X$ be a smooth algebraic surface over $k$.
We recall the follwoing notations from section~\ref{surface}.
If $C$ is a curve on $X$, then $K_C$ is the completion of the field
$k(X)$ along the discrete valuation given by the irreducible curve $C$ on $X$.
If $p$ is a point on $X$, then $K_p$ is a subring in
$ \Frac(\hat{\oo}_{p, X}  )$ generated by subrings $k(X)$ and $\hat{\oo}_{p,X}$.

There are the following reciprocity laws, see~\cite{FP}.
\begin{Th}
\begin{enumerate}
\item We fix a point $p \in X$ and any $f, g, h \in K_p^*$.
Then in the following sum only a finite number of terms is
not equal to $0$ and
$$
\sum_{ X \supset \C \ni p} \nu_{K_{p,\C}} (f,g) = 0 \mbox{.}
$$
In the following product only a finite number of terms is not equal to $1$ and
$$
\prod_{X \supset \C \ni p} (f,g,h)_{K_{p, \C}} = 1 \mbox{,}
$$
where the sum and the product are taken over all germs of irredusible
curves on $X$ at $p$.
\item
We fix an irreducible projective curve $C$ on $X$ and any $f,g,h \in K_C^*$.
Then in the following sum only a finite number of terms is
not equal to $0$ and
$$
\sum_{   p \in \C \subset X } [k(p) : k] \cdot \nu_{K_{p,\C}} (f,g) = 0 \mbox{.}
$$
In the following product only a finite number of terms is not equal to $1$ and
$$
\prod_{ p \in \C \subset X} \nm\nolimits_{k(p) /k} \, (f,g,h)_{K_{p, \C}} = 1 \mbox{,}
$$
where the sum and the product are taken over all points $p \in C$  and all germs of  irredusible
curve $C$ at $p$.
\end{enumerate}
\end{Th}

\begin{nt} {\em
The relative reciprocity laws were constructed in~\cite{Os0} (see also \cite{Os00} for a
short exposition)
for a smooth projective morphism $f$ of a smooth algebraic surface $X$ to
a smooth algebraic curve $S$ when $char k = 0$.
 If $p \in \C \subset X$ then explicit formulas were constructed in~\cite{Os0} for
maps
$$
K_2 (K_{p, \C}) \lto K_{f(p)}^*  \mbox{.}
$$
 }
\end{nt}

\begin{nt} {\em
The map $\nu_K (\;, \;)$ for $2$-dimensional local field $K$
was interpreted in~\cite{Os1} as the commutator of liftings of elements $f, g \in K^*$
in a central extension of group $K^*$ by $\dz$. From this interpretation
the reciprocity laws for $\nu_K (\; , \;)$ were proved.
The proof in~\cite{Os1} uses  adelic rings on an algebraic surface $X$.
This is an abstract version of reciprocity law for $\nu_K(\; , \;)$ like abstract
version of reciprocity law~(\ref{curver}) for a projective curve  in~\cite{Arbar} (and in~\cite{T} for the residues of differentials on a projective curve).}
\end{nt}

\begin{nt} {\em
We don't describe
the reciprocity laws for residues of differentials of $2$-dimensional local fields, which were
formulated and proved in~\cite{P1}, see also~\cite{Y}. }
\end{nt}

\begin{nt} {\em
The symbols $\nu_K (\; , \;)$ and $(\;, \;, \;)_K$
correspond to the non-ramified and tame ramified extensions of $2$-dimensional local fields
when the last residue field is finite, see~\cite{P}.
}
\end{nt}

\noindent Steklov Mathematical Institute, \\
Gubkina str. 8, \\
119991, Moscow, Russia \\
e-mail  ${ d}_{-}osipov@mi.ras.ru$

\end{document}